\documentclass[a4]{article}

\usepackage{amsmath,amsfonts,amsthm,amssymb,graphicx,tikz-cd,stmaryrd,enumitem}
\allowdisplaybreaks

\usetikzlibrary{graphs}
\theoremstyle{definition}
\newtheorem{thm}{Theorem}[section]
\newtheorem{deff}[thm]{Definition}
\newtheorem{lemm}[thm]{Lemma}
\newtheorem{prop}[thm]{Proposition}
\newtheorem{cor}[thm]{Corollary}
\newtheorem{rem}[thm]{Remark}
\newtheorem{exam}[thm]{Example}

\newcommand{\zahl}{\mathbb{Z}}

\newcommand{\card}{\sharp}

\newcommand{\pone}{\mathbb{P}^1}

\newcommand{\ol}[1]{\overline{#1}}
\newcommand{\fct}[1]{\operatorname{#1}}

\newcommand{\slt}{\fct{SL}_2}

\newcommand{\eqn}[1]{\[ #1 \]}
\newcommand{\eqnl}[1]{\begin{equation} #1 \end{equation}}
\newcommand{\eqns}[1]{\begin{align*} #1 \end{align*}}
\newcommand{\eqnsl}[1]{\begin{align} #1 \end{align}}

\title{Bracket Polynomial Expression of Discriminant-Resultants as $\slt$-invariant}
\author{Rin Gotou\thanks{Department of Mathematics, Graduate School of Osaka University \texttt{u661233h@ecs.osaka-u.ac.jp}}}
\date{}

\begin{document}
\maketitle
\begin{abstract}
We give a bracket polynomial expression for intermediate terms between discriminant and resultant for pair of binary forms. As an application of the bracket polynomial expression, we give an algebraic proof of the algebraic independence of intermediate terms, which was shown in the theory of dynamical systems.
\end{abstract}
\section{Introduction}
Throughout this paper, we fix a field $k$ of characteristic zero. For any nonnegative integer $m$, let $V_m := k[x,y]_m$ be the vector space of the binary forms of degree $m$. For a vector space $V$, we denote by $k[V]$ the algebra $\fct{Sym}^\bullet(V^*)$ of the regular functions over $V$.

% For a given $\slt$-representation $V$, we consider $k[V]^{\slt}$. There are several methods to construct invariants from covariants, elements in $k[V \oplus V_1]^{\slt}$. The self-duality of the representation $V_1$ gives the \emph{fundamental covariant} of $V_n$, the unique covariant $f_n \in k[V_n \oplus V_1]_{(1,n)} \simeq V_n \otimes V_n^*$ up to constant multiple.

%The isomorphism between dual representations $V_1^* \to V_1$ by $x^* \mapsto y,\ y^* \mapsto -x$ gives the \emph{fundamental covariant} of $V_n$,
% \[f_n(x,y) := \sum_{i=0}^n a_i \otimes \binom{n}{i}(-1)^{n-i}x^{n-i} y^i \in V_n \otimes \fct{Sym}_{n}V_1 . \]
Discriminants and resultants are classical invariants which play a fundamental role in algebraic geometry (\cite{C-L-O},\cite{C-L-O2}). Here we note some elemental properties. Let $d$ and $e$ be any positive integers. The resultant $\fct{res} = \fct{res}_{d,e}$ (for binary forms of degree $d$ and $e$) is a vanishing polynomial of the locus $\{ (f,g) \in V_d \times V_e \mid f\text{ and }g \text{ have a common factor} \}$ in the $V_d \times V_e$ as an affine space. Here the variables consists the polynomial are elements in $V_d^*$ and $V_e^*$. We have $\fct{res} \in k[V_d \times  V_e]_{(e,d)}$ with the natural multi-grading on $k[V_d \times V_e]$  given by the elements in $V_d^*$ and $V_e^*$ to be $(1,0)$ and $(0,1)$ respectively. The discriminant is a vanishing  polynomial of the locus $\{ f \in V_d \mid f \text{ has a multiple factor} \}$. The discriminant $\Delta \in k[V_{d}]_{2d-2}$ is given by $\Delta(f) := \fct{res}(f,x\partial_{x}f)/a_0a_d\ (\text{where }f = \sum_{i = 0}^d a_ix^iy^{d-i})$. Under the $\slt$-action on $V_d \times V_e$ as the $\slt$-representation $V_d \oplus V_e$, the resultant and the discriminant are $\slt$-invariant.

In \cite{G}, the author found that intermediate terms between discriminant and resultant plays a fundamental role in representation-theoretic treatment of the moduli space of dynamical systems over the project line. In this paper, we call the intermediate terms \emph{discriminant-resultant}. Let $n \geq 2$ be a positive integer. 
\begin{deff}
The $r$-th discriminant-resultants $DR_{n,r} \in k[V_n \times V_{n-2}]_{(2n-2-r,r)}$ for $0 \leq r \leq n$ is the polynomials which satisfy
\[ \sum_{r=0}^n DR_{n,r}(f_n,f_{n-2})t^r = \fct{res}(f_n(x,y), x\partial_xf_n(x,y) + txy f_{n-2}(x,y)) / a_0a_n. \]
\end{deff}
It is easy to see $DR_{n,0} = \Delta(f_n)$ and $DR_{n,n} = \fct{res}(f_n , f_{n-2})$.
In \cite{G}, the author showed that $DR_{n,1} = 0$ and $DR_{n,r} (2 \leq r \leq n-1)$ are $\slt$-invariant. Despite discriminant-resultants are $\slt$-invariant, the theory of generalizations of resultants are considered mainly under the toric action on the variables ($\mathcal{A}$-resultant in \cite{G-K-Z},\cite{C-L-O}).

\begin{rem} The discriminant-resultants are a special case of van der Wearden's $u$-resultant (\cite{C-L-O}). The simplest type of $u$-resultant is defined like $\fct{res}(f , gt + h)$. The main purpose of the theory of $u$-resultant is to complete the universality of the theory of resultant. Argument about this special case of $u$-resultant is less required before the relation between dynamical system have been found.
\end{rem}
To treat discriminant-resultants as $\slt$-invariants, we use \emph{the symbolic method} \cite[Chapter 9]{D-C}, a classical universal process to construct $\fct{SL}_N$-invariants. Let $V \simeq \bigoplus_j V_{m_j , \alpha_j}$ be an $\slt$-representation, where $V_{m_j,\alpha_j}$ is $\slt$-representation isomorphic to $V_{m_j}$, indexed by a symbol $\alpha_j$ to distinguish the summands. We write $f_{m,\alpha}(x,y)$ for a generic element of $V_{m,\alpha}$. The symbolic method is done by introducing $(m \times 2)$-tuple of variables $\alpha_{i,\epsilon} (1 \leq i \leq m, \epsilon = 0,1)$ for each representation $V_{m,\alpha}$, which assumed to satisfy the following equation for the binary form:
\begin{equation}
 f_{m,\alpha}(x,y) = \sum_{i = 0}^m a_ix^iy^{m-i} = \prod_{i=1}^m(\alpha_{i,0} x - \alpha_{i,1} y ) \label{alpha}.
\end{equation}
For any indices of representations $\alpha$ and $\beta$ (possibly $\alpha = \beta$), \emph{the bracket symbols} \[ [\alpha_i , \beta_j] := \alpha_{i,0}\beta_{j,1} - \beta_{j,0}\alpha_{i,1} \]
are virtual $\slt$-invariants. Any finite product of bracket symbols is called \emph{bracket monomial}. Some symmetric sums of bracket monomials belong to $k[V]^{\slt}$ under the relation (\ref{alpha}). Such symmetric sums can express all elements in $k[V]^{\slt}$ (see Kung and Rota (\cite{K-R}) for a proof). 
\begin{exam}
On the representation $V = V_{d,\alpha} \oplus V_{e,\beta}$, the resultant is known to be
\[ \fct{res}(f_{d,\alpha}(x,y),f_{e,\beta}(x,y)) = \prod_{i = 1}^d \prod_{j = 1}^e [\alpha_i ,\beta_j ]  \]
and the discriminant is known to be
\[ \Delta (f_{d,\alpha}(x,y)) = \prod_{1 \leq i,j \leq d, i \neq j} [\alpha_i,\alpha_j]. \]
\end{exam}
% \begin{rem}
% From the viewpoint of the geometric invariant theory, the symbolic method gives rational functions over $\proj (\bigoplus_{i} V_{n_i})(\oshf (1) )^{ss} \sslash \slt$ through the isomorphism
% \eqns{\proj \left( \bigoplus_{i} V_{n_i} \right)(\oshf (1) )^{ss} \sslash \slt & \simeq \proj \left( V_1^{\oplus \sum_i n_i} \right)(\oshf(1,\ldots , 1))^{ss} \sslash \prod_i \sym_{n_i} \times \slt  \\
% &\simeq G\left( \sum_{i}n_i , 2\right)((Pl)^* \oshf(1))^{ss} \sslash \prod_i \sym_{n_i},}
% where $Pl$ is the Pl\"{u}cker embedding.
% \end{rem}

In this paper, we give the following bracket polynomial expression of discriminant-resultants. 
\begin{thm}(Theorem \ref{CayleyDR})\label{main1}
    We have $DR_{2,2} = f_0^2$ and
    \begin{equation}
         DR_{n,r} = \sum_{\substack{ I \sqcup J = [n], \\|I| = r } }\left( \prod_{\substack{ j \in J \\ i \in [n] \setminus \{ j \}}} [\alpha_i , \alpha_j ] \cdot \prod_{\substack{i \in I \\ k \in [n-2]} } [\beta_k , \alpha_i] \right)
         \label{eqmain1}
    \end{equation}
    except $(n,r) = (2,2)$, where $[m] := \{ 1,2,\ldots , m\}$. 
\end{thm}

As an application of Theorem \ref{main1}, we show two results about discriminant-resultants by purely algebraic method. These results are already obtained from the theory of dynamical systems in \cite{G}. One is the vanishing of the first discriminant-resultant.
\begin{thm} (Corollary \ref{drnone}) \label{main2} We have $DR_{n,1} = 0$.
\end{thm}
\begin{rem}The author commented two other proofs of Theorem \ref{main2} in \cite{G}. One is obtained as a representational interpretation of the Woods Hole formula \cite{A-B}. The other is directly obtained from an interpretation of the definition as a deformation of the discriminant.
\end{rem}
The other is the algebraic independence of discriminant-resultants, which was shown as the non-degeneracy of fixed point multiplier map of dynamical systems over the projective line. The case of $n = 3$ is shown by Milnor (\cite{Mi}) and $n \geq 4$ is shown by Fujimura (\cite{Fu2}). In \cite{G}, the author translated their result into the terms of discriminant-resultant.
\begin{thm}(Theorem \ref{algindep})\label{main3}
The discriminant-resultants \[ \{ DR_{n,r} \mid r=0,2,3,\ldots ,n \} \] are algebraically independent.
\end{thm}
To show Theorem \ref{algindep}, we use some elemental facts which motivates the theory of cluster algebra (\cite{F-Z}). The theory of cluster algebra originates in an observation that the Grassmannian algebra, the algebra of bracket symbols, have various embeddings into Laurent rings. Despite cluster algebra is constructed from such many ways of embeddings and relations among them, we only use one fixed embedding. Then, we take a monomial order on the Laurent ring and show Theorem \ref{algindep} using the monomial order.

This paper is organized as the following. In Section 2, we give a proof of Theorem \ref{main1}. In Section 3, we define the Grassmannian algebra, the algebra of Cayley symbols. We also introduce some elemental fact about embeddings of the Grassmannian algebra into Laurent rings. In Section 4, we show Theorem \ref{main2}. In Section 5, we introduce a monomial ordering on a Laurent ring and we show Theorem \ref{main3}.
\subsection*{Acknowledgments}
The author was supported by JSPS KAKENHI Grant-in-Aid for Research Fellow JP202122197. 

\section{Proof of Theorem \ref{main1}}
\begin{thm}\label{CayleyDR}
    We have $DR_{2,2} = f_0^2$ and
    \begin{equation}
         DR_{n,r} = \sum_{\substack{ I \sqcup J = [n], \\|I| = r } }\left( \prod_{\substack{ j \in J \\ i \in [n] \setminus \{ j \}}} [\alpha_i , \alpha_j ] \cdot \prod_{\substack{i \in I \\ k \in [n-2]} } [\beta_k , \alpha_i] \right)
         \label{eqCDR}
    \end{equation}
    except $(n,r) = (2,2)$, where $[m] := \{ 1,2,\ldots , m\}$. 
\end{thm}
\begin{proof}

We show the assertion by an induction on $n$. The base cases are the cases $r = 0$ of discriminants, the cases $r=n$ of resultants and the cases $(n,r) = (2,1),(3,2)$. The cases $r=0$ and $r=n$ are already introduced. The cases $(n,r) = (2,1),(3,2)$ are shown by the direct computations.

We assume that (\ref{eqCDR}) holds true for $DR_{n-1,r}$ for any $r$. Firstly, we consider the case that $f_n$ and $f_{n-2}$ have a common linear factor $h = (px + qy)$. In this case, the assertion is obtained by direct computations on the both hand of the equation (\ref{eqCDR}), as the following. We put $f_{n-1}$ and $f_{n-3}$ as $f_n = f_{n-1} \cdot h$ and $f_{n-2} = f_{n-3} \cdot h$ respectively.
Then, properties of resultant admits us to compute the resultant with $f_n$ and any polynomial $F$ as
\eqns{
   \fct{res} ( f_n , F)
   &=\fct{res}(f_{n-1} h, F) \\
   &=\fct{res}(f_{n-1}, F)\cdot \fct{res}(h, F) \\
   &=\fct{res}(f_{n-1}, F\ \fct{mod}\ f_{n-1})\cdot \fct{res}(h, F\ \fct{mod}\ h).}
We apply this for $F = x \partial_x f_n + txyf_{n-2}$. Here, we have 
\eqns{
x \partial_x f_n + txyf_{n-2} &= x\partial_x (hf_{n-1}) + txyhf_{n-3} \\
&= x((\partial_x h) \cdot f_{n-1} +  h (\partial_xf_{n-1})) + txyhf_{n-3}\\
&= \begin{cases}
h( x\partial_xf_{n-1} + txyf_{n-3} ) & \mod f_{n-1} \\
x(\partial_x h) \cdot f_{n-1} & \mod h. 
\end{cases}
}
Therefore, we obtain
\eqns{
   \fct{res} ( f_n , x \partial_x f_n + txyf_{n-2} ) 
=&  \fct{res}(f_{n-1},h( x\partial_xf_{n-1} + txyf_{n-3} )) \cdot \fct{res}( h , x (\partial_x h) \cdot f_{n-1}) \\
=& \fct{res}(h,x (\partial_x h) ) \cdot \fct{res}( h, f_{n-1}) \cdot \fct{res}(f_{n-1} , h)  \\ 
   & \cdot  \fct{res}(f_{n-1},x\partial_x f_{n-1} + txyf_{n-3}) \\
=& (-1)^{n-1}pq \cdot \fct{res}( h, f_{n-1})^2 \cdot  \fct{res}(f_{n-1},x\partial_x f_{n-1} + txyf_{n-3}).
}
We write $f_n = \sum_{i = 0}^n a_ix^iy^{n-i}$ and $f_{n-1} = \sum_{i=0}^{n-1} a'_ix^iy^{n-1-i}$. Then $f_n = f_{n-1} \cdot h$ leads to $a_0= pa'_0$ and $a_n = qa'_{n-1}$. Therefore, we have 
\eqns{ \sum_{r=0}^n DR_{n,r} t^r &
= \fct{res} ( f_n , x \partial_x f_n + txyf_{n-2} )/(a_0a_n) \\
& = (-1)^{n-1} pq \cdot \fct{res}( h, f_{n-1})^2 \cdot  \fct{res}(f_{n-1},x\partial_x f_{n-1} + txyf_{n-3})/(a_0a_n) \\
& =(-1)^{n-1}\fct{res}( h, f_{n-1})^2 \cdot \fct{res}(f_{n-1},x\partial_x f_{n-1} + txyf_{n-3})/(a'_0a'_{n-1}) \\
& =(-1)^{n-1}\fct{res}( h, f_{n-1})^2\cdot \sum_{r=0}^{n-1} DR_{n-1,r}t^r.}
Therefore, we obtain
\eqn{ DR_{n,r} = (-1)^{n-1}\fct{res}( h, f_{n-1})^2 \cdot DR_{n-1,r}.}
On the RHS of (\ref{eqCDR}), we have $\alpha_n = \beta_{n-2} = h$ and $[\alpha_n , \beta_{n-2}] = 0$ from the assumption. 
Thus, we have
\eqnl{
(\text{the RHS of }(\ref{eqCDR})) = \sum_{\substack{ I \sqcup J = [n], \\ \card{I} = r \\ n \not \in I} }\left( \prod_{\substack{ j \in J \\ i \in [n] \setminus \{ j \}}} [\alpha_i , \alpha_j ] \cdot \prod_{\substack{i \in I \\ k \in [n-2]} } [\beta_k , \alpha_i] \right). \label{rhseqn}}
For each term in (\ref{rhseqn}), we have $n \in J$ and thus
\eqns{
\prod_{\substack{ j \in J \\ i \in [n] \setminus \{ j \}}} [\alpha_i , \alpha_j ] & =  \prod_{\substack{ j \in J' \\ i \in [n-1] \setminus \{ j \}}} [\alpha_i , \alpha_j ]  \cdot \prod_{j \in J'} [\alpha_n,\alpha_j] \cdot \prod_{i \in [n-1]} [\alpha_i,\alpha_n]  \text{ and }\\
\prod_{\substack{i \in I \\ k \in [n-2]} } [\beta_k , \alpha_i] & = \prod_{\substack{i \in I \\ k \in [n-3]} } [\beta_k , \alpha_i] \cdot \prod_{ i \in I } [\beta_{n-2} , \alpha_i].
}
Here we have
\eqn{
\prod_{j \in J'} [\alpha_n,\alpha_j] \cdot \prod_{i \in [n-1]} [\alpha_i,\alpha_n] \cdot \prod_{ i \in I } [\beta_{n-2} , \alpha_i] = (-1)^{n-1} \left( \prod_{i \in [n-1]}[\alpha_i,h]^2 \right)
}
in any term, thus 
\eqns{
(\ref{rhseqn})& = (-1)^{n-1} \left( \prod_{i \in [n-1]}[\alpha_i,h]^2 \right) \cdot \sum_{\substack{ I \sqcup J' = [n-1], \\ \card{I} = r}}\left( \prod_{\substack{ j \in J' \\ i \in [n] \setminus \{ j \}}} [\alpha_i , \alpha_j ] \cdot   \prod_{\substack{i \in I \\ k \in [n-3]} } [\beta_k , \alpha_i] \right) \\
& = (-1)^{n-1} \fct{res}(f_{n-1} , h)^2 \cdot DR_{n-1,r}.
}
Therefore, the assertion is true if $f_n$ and $f_{n-2}$ have a common linear factor. Since the pair $f_n$ and $f_{n-2}$ have a common linear factor if and only if $\fct{res}(f_n , f_{n-2}) = 0$, and the fact that resultant is an irreducible polynomial (see \cite{G-K-Z}), we have
\eqnl{  DR_{n,r} - \sum_{\substack{ I \sqcup J = [n], \\|I| = r } }\left( \prod_{\substack{ j \in J \\ i \in [n] \setminus \{ j \}}} [\alpha_i , \alpha_j ] \cdot \prod_{\substack{i \in I \\ k \in [n-2]} } [\beta_k , \alpha_i] \right) = g \cdot \fct{res}(f_n, f_{n-2})
\label{remres}}
for some polynomial $g$. As a polynomial of coefficients of $f_{n-2}$, all terms on the LHS of (\ref{remres}) is homogeneous of degree $r$ and $\fct{res}(f_n,f_{n-2})$ is of degree $n$. Therefore, we have $g=0$ for $r \leq n-1$. This completes the induction.
\end{proof}

\section{Laurent Phenomenon for Cluster Algebra}
For a finite set of symbols $S$, we write $G(S)$ for the algebra generated by the brackets of symbols. The relations among the bracket symbols are generated by the anti-commutation relations $[\alpha , \beta ] = -[\beta , \alpha ]$ and the Pl\"{u}cker relations
\[ [\alpha , \beta] [\gamma , \delta ] + [\alpha , \gamma] [\delta , \beta ] +[\alpha , \delta ] [\beta , \gamma ] = 0 \]
for any symbols $\alpha , \beta , \gamma $, and $ \delta$.
Therefore, the algebra $G(S)$ is isomorphic to the Grassmannian algebra $k [ G(2, |S|) ]$, the ring generated by the Pl\"{u}cker coordinates of the Grassmannian variety $G(2,|S|)$.

The Grassmannian algebras $k [ G(2,n) ]$ are motivating examples of the theory of cluster algebra (see \cite{F-Z} and \cite{Sc}).
To show the algebraic independence of discriminant-resultants, we use valuations induced from the ``Laurent Phenomena'' of the Grassmannian algebra.

To explain the Laurent Phenomena and their conclusions, it is convenient to place the symbols $S$ on a vertices of a convex $\card{S}$-gon $P$ on plane. We associate the bracket symbol $[\alpha , \beta ]$ with the line segment $\ol{\alpha \beta}$. Once we fix a diagonal $\ol{\alpha\beta}$ of the polygon $P$, we can write any brackets by a Laurent polynomial of other brackets whose line segment does not crosses $\ol{\alpha \beta}$ at the interior of the polygon $P$, because if $\ol{\gamma \delta}$ crosses $\ol{\alpha \beta}$ then we have 
\eqnl{ [\gamma , \delta ]= -\frac{ [\alpha , \gamma] [\delta , \beta ] +[\alpha , \delta ] [\beta , \gamma ]}{[\alpha , \beta]}. \label{plured} }
By repeating this process for the diagonals given by a fixed triangulation $T$ of the polygon $P$, we obtain the following.
%Then, for symbols $\alpha , \beta , \gamma $ and $\delta$ placed counter-clockwise on the vertices, the Pl\"{u}cker relation
% \[ [\alpha , \gamma] [\beta , \delta] =  [\alpha , \beta] [\gamma , \delta ] +[\alpha , \delta ] [\beta , \gamma ] \]
% works as resolving the intersection of segments $\ol{\alpha \gamma}$ and $\ol{\beta \delta}$. Moreover, if we fix a diagonal $\ol{\alpha\gamma}$ of the polygon, we can deform any bracket  crossing $\ol{\alpha \gamma}$
% \[ [\beta , \delta] = \left( [\alpha , \beta] [\gamma , \delta ] +[\alpha , \delta ] [\beta , \gamma ] \right) / [\alpha , \gamma ]. \]

\begin{lemm}(Laurent Phenomena of Cluster Algebra)\label{laulemm} Let $S$ be a finite set of symbols and assume that $\card{S} \geq 3$. Let $P$ be a convex $n$-gon on a plane with vertices indexed by the elements of $S$. Let $T$ be a triangulation of $P$, $D_T$ the $n-3$ diagonals of $P$ defining the triangulation $T$ and $E_T$ the $2n-3$ line segments of triangles of $T$, that is, $E_T = \{ \text{edges of } P\} \cup D_T$. Then, we have the following:
\begin{enumerate}[label = (\roman*)]
    \item \label{lauphen} The Pl\"{u}cker relations leads the natural inclusion \begin{equation}
        k[ E_T ] \subset G(S) \subset k [E_T,D_T^{-1}].
    \end{equation} 
    \item \label{laurep} Let $\alpha = \delta_0 , \delta_1 , \ldots , \beta = \delta_k $ and $\gamma$ be distinct symbols in $S$. Then we have
    \[ [\alpha , \beta ] = \sum_{i=0}^{k-1} \frac{ [\gamma , \alpha ][ \gamma, \beta ][\delta_i , \delta_{i+1} ] }{ [\gamma , \delta_i ] [\gamma, \delta_{i+1} ] }. \]
\end{enumerate}
\end{lemm}
\begin{proof} \ref{lauphen} See \cite{F-Z}.
\ref{laurep} This is obvious from (\ref{plured}) and $[\alpha , \beta] = -[\beta , \alpha ]$.
\end{proof}
% For a pair $(\gamma , \delta )$ of symbol and 
% \begin{lemm}

% \end{lemm}
\section{Proof of Theorem \ref{main2}}
We show Theorem \ref{main2} by a direct computation using Theorem \ref{main1}.
\begin{cor} We have $DR_{n,1} = 0$.
\end{cor}
\begin{proof}
We show the assertion by the induction on $n$. By a direct computation, we have $DR_{2,1} = 0$. We assume $DR_{n-1,r} = 0$. By Theorem \ref{main1}, we have
\eqn{ \sum_{l = 1}^{n-1} \prod_{\substack{ j \in [n-1] \setminus \{ l \} \\ i \in [n-1] \setminus \{ j \}}}[\alpha_i , \alpha_j ] \cdot \prod_{k \in [n-3]}[\beta_k , \alpha_l] = 0. }
We transpose this as
\eqnl{ \prod_{\substack{ j \in [n-1] \setminus \{ n-1 \} \\ i \in [n-1] \setminus \{ j \}}}[\alpha_i , \alpha_j ] \cdot \prod_{k \in [n-3]}[\beta_k , \alpha_{n-1}] = 
-\sum_{l = 1}^{n-2} \prod_{\substack{ j \in [n-1] \setminus \{ l \} \\ i \in [n-1] \setminus \{ j \}}}[\alpha_i , \alpha_j ] \cdot \prod_{k \in [n-3]}[\beta_k , \alpha_l] . \label{tped} }
Again by Theorem \ref{main1}, we have
\eqnl{ DR_{n,1} = \sum_{l = 1}^n \prod_{\substack{ j \in [n] \setminus \{ l \} \\ i \in [n] \setminus \{ j \}}}[\alpha_i , \alpha_j ] \cdot \prod_{k \in [n-2]}[\beta_k , \alpha_l]. \label{drnone} }
We note that, for any $1 \leq l \leq n-1$,
\eqns{\prod_{\substack{ j \in [n] \setminus \{ l \} \\ i \in [n] \setminus \{ j \}}}[\alpha_i , \alpha_j ] &= \prod_{\substack{ j \in [n-1] \setminus \{ l \} \\ i \in [n-1] \setminus \{ j \}}}[\alpha_i , \alpha_j ] \cdot \prod_{\substack{ i \in [n-1] }}[\alpha_i , \alpha_n ] \cdot \prod_{\substack{ j \in [n-1] \setminus \{ l \} }}[\alpha_n , \alpha_j ]
\\ &= \frac{\prod_{\substack{ i \in [n-1] }}\left( -[\alpha_n,\alpha_i]^2 \right) }{[\alpha_n,\alpha_l]}. \cdot  \prod_{\substack{ j \in [n-1] \setminus \{ l \} \\ i \in [n-1] \setminus \{ j \}}}[\alpha_i , \alpha_j ]. } 
 In particular, the term of $l=n-1$ in the RHS of (\ref{drnone}) have the LHS of (\ref{tped}) as a factor. By substituting it, we obtain
\eqnsl{
& \prod_{\substack{ j \in [n] \setminus \{ n-1 \} \\ i \in [n] \setminus \{ j \}}}[\alpha_i , \alpha_j ] \cdot \prod_{k \in [n-2]}[\beta_k , \alpha_{n-1}] \nonumber \\
& =\frac{\prod_{\substack{ i \in [n-1] }}\left( -[\alpha_n,\alpha_i]^2 \right) }{[\alpha_n,\alpha_{n-1}]} \cdot \prod_{\substack{ j \in [n-1] \setminus \{ n-1 \} \\ i \in [n-1] \setminus \{ j \}}}[\alpha_i , \alpha_j ] 
\cdot \prod_{k \in [n-2]}[\beta_k , \alpha_{n-1}] \nonumber \\
&= \frac{\prod_{\substack{ i \in [n-1] }}\left( -[\alpha_n,\alpha_i]^2 \right) }{[\alpha_n,\alpha_{n-1}]}
\cdot \left( -\sum_{l = 1}^{n-2} \prod_{\substack{ j \in [n-1] \setminus \{ l \} \\ i \in [n-1] \setminus \{ j \}}}[\alpha_i , \alpha_j ] \cdot \prod_{k \in [n-3]}[\beta_k , \alpha_l] \right) \cdot [\beta_{n-2},\alpha_{n-1}] \nonumber \\
&= \left( -\sum_{l = 1}^{n-2} \frac{[\alpha_n , \alpha_l] }{[\alpha_n,\alpha_{n-1}]} \prod_{\substack{ j \in [n] \setminus \{ l \} \\ i \in [n] \setminus \{ j \}}}[\alpha_i , \alpha_j ] \cdot \prod_{k \in [n-3]}[\beta_k , \alpha_l] \right) \cdot [\beta_{n-2},\alpha_{n-1}]. \label{termone}
}
By the similar substitution for the term of $l = n$ in the RHS of (\ref{drnone}), we obtain 
\eqnsl{
& \prod_{\substack{ j \in [n] \setminus \{ n \} \\ i \in [n] \setminus \{ j \}}}[\alpha_i , \alpha_j ] \cdot \prod_{k \in [n-2]}[\beta_k , \alpha_{n-1}] \nonumber \\
&= \left( -\sum_{l = 1}^{n-2} \frac{[\alpha_{n-1} , \alpha_l] }{[\alpha_{n-1},\alpha_n]} \prod_{\substack{ j \in [n] \setminus \{ l \} \\ i \in [n] \setminus \{ j \}}}[\alpha_i , \alpha_j ] \cdot \prod_{k \in [n-3]}[\beta_k , \alpha_l] \right) \cdot [\beta_{n-2},\alpha_n]. \label{termtwo}
}
By the Pl\"{u}cker relation (\ref{plured}), we have 
\[ \frac{[\alpha_n , \alpha_l] [\beta_{n-2},\alpha_{n-1}] }{[\alpha_n,\alpha_{n-1}]} + \frac{[\alpha_{n-1} , \alpha_l] [\beta_{n-2},\alpha_n] }{[\alpha_{n-1},\alpha_n]} = [\beta_{n-2},\alpha_l] \]
for each $1\leq l \leq n-2$. Therefore, we have
\[ (\ref{termone}) + (\ref{termtwo}) = - \sum_{l = 1}^{n-2} \prod_{\substack{ j \in [n] \setminus \{ l \} \\ i \in [n] \setminus \{ j \}}}[\alpha_i , \alpha_j ] \cdot \prod_{k \in [n-3]}[\beta_k , \alpha_l] \cdot [\beta_{n-2},\alpha_l], \]
this shows that $DR_{n,1} = 0$.
\end{proof}

\section{Proof of Theorem \ref{algindep}}
To show Theorem \ref{algindep}, we use the polygon whose vertices are marked by the symbols $\alpha_1 , \ldots , \alpha_n , \beta_1 , \ldots , \beta_{n-2}$ counter-clockwise, and the triangulation $T$ given by drawing all diagonals through $\beta_{n-2}$. We put
\begin{align*}
    A_i & := [\beta_{n-2},\alpha_{i}] & & (i = 1, \ldots , n) \\
    B_i & := [\beta_{n-2},\beta_{i}] & & (i = 1, \ldots , n-3)\\
    C_i & := \begin{cases}[\alpha_i,\alpha_{i+1}] \\
    [\alpha_n,\beta_1]
    \end{cases}& & \begin{aligned}
         & (i = 1 , \ldots ,  n-1) \\
         & (i = n)  
    \end{aligned}\\
    D_i & := [\beta_i, \beta_{i+1}] & & (i = 1, \ldots , n-4).
\end{align*}
The Laurent polynomial ring corresponding to the triangulation $T$ is 
\[ L := k [A_i , B_i , C_i, D_i , A_2^{-1},\ldots , A_n^{-1},B_1^{-1},\ldots, B_{n-4}^{-1}].\]
By Lemma \ref{laulemm} \ref{lauphen}, we regard $G(\{ \alpha_i , \beta_i \} )$ as a subring of $L$.

Our strategy is to fix a monomial ordering on $L$ and see the leading monomials of $DR_{n,r}$. For a Laurent polynomial $f \in L$ and a monomial ordering $\preceq$ on $L$, we write the leading monomial of $f$ with respect to the ordering $\preceq$ by $\fct{lm}_{\preceq}(f)$.

\begin{prop}\label{lm} Let $\preceq$ be a monomial ordering and $\{ f_s \}_{s \in S}$ be a finite family of Laurent polynomials in $L$. If there is a polynomial $f_s$ such that $\fct{lm}_{\preceq}(f_s) \succ \fct{lm}_{\preceq}(f_{s'})$ for any other $s' \in S$, then we have $\fct{lm}_{\preceq}(\sum_{s \in S} f_s) = \fct{lm}_{\preceq}(f_s)$.
\end{prop}
\begin{deff}A subset $S$ of $L$ is said to be \emph{multiplicatively independent} if the map \[ \bigoplus_{s \in S}\zahl \ni (a_s)_{s \in S} \mapsto \prod_{\substack{s \in S \\ a_s \neq 0}} s^{a_s} \in L \] is injective.
\end{deff}
\begin{cor}\label{mulindep} A family of Laurent polynomials $\{ f_s \}_{s \in S}$ are algebraically independent if there exists a monomial ordering $\preceq$ on $L$ such that $\fct{lm}_{\preceq}(f_s)$ are multiplicatively independent.
\end{cor}
We fix the monomial ordering $\preceq$ on $L$ as the lexicographic order, where the variables are sorted as \[ ( A_1,\ldots , A_n,B_1,\ldots , B_{n-3}, C_1,\ldots , C_n , D_1,\ldots , D_{n-4}) \]
from maximal to minimal.
To show Theorem \ref{algindep}, we see the degree of $A_i$'s and $C_i$'s on $\fct{lm}_{\preceq}(DR_{n,r})$.
\begin{lemm}
We have 
\eqns{ \fct{lm}_{\preceq} ([\alpha_i , \alpha_j ]) &= A_i A_{j-1}^{-1}C_{j-1} \ (i < j) \text{ and}
\\ \fct{lm}_{\preceq} ([\alpha_i , \beta_k ]) &= \begin{cases} A_iA_n^{-1}C_n & (k = 1) \\ A_iB_{k-1}^{-1}D_{k-1}& (2\leq k \leq n-3) \\
A_i & (k = n-2). \end{cases}\\}
\end{lemm}
\begin{proof}
The assertion is obvious from Lemma \ref{laulemm} \ref{laurep}.
\end{proof}
\begin{lemm}\label{mondegAl}
Let $l \in [n]$ and $I\subset [n]$. We put $J := [n] \setminus I$ and $r := \card I $. Then, we have 
\eqns{
& \deg_{A_l} \fct{lm}_{\preceq} \left(  \prod_{\substack{ j \in J \\ i \in [n] \setminus \{ j \}}} [\alpha_i , \alpha_j ] \cdot \prod_{\substack{i \in I \\ k \in [n-2]} } [\beta_k , \alpha_i] \right)\\
& = \begin{cases}n-2 & (l \in I) \\ n-l & (l \in J) \end{cases} - \begin{cases} 0 & (l+1 \in I) \\ l & (l+1 \in J) \end{cases} + \begin{cases}n-r - 2\cdot \card{ ( J \cap [l] ) } & (l \leq n-1) \\ -r & (l = n).\end{cases}
}
\end{lemm}
\begin{proof} For convenience, we put $B_0 := A_n$ and $D_0 := C_n$. Then we have
\eqnsl{
& \fct{lm}_{\preceq} \left(  \prod_{\substack{ j \in J \\ i \in [n] \setminus \{ j \}}} [\alpha_i , \alpha_j ] \cdot \prod_{\substack{i \in I \\ k \in [n-2]} } [\beta_k , \alpha_i] \right) \nonumber \\
& = \prod_{\substack{ j \in J \\ i \in [n] \setminus \{ j \}}} \fct{lm}_{\preceq}([\alpha_i , \alpha_j ]) \cdot \prod_{\substack{i \in I \\ k \in [n-2]} } \fct{lm}_{\preceq} ( [\beta_k , \alpha_i] )  \nonumber \\
& = \prod_{\substack{ j \in J}} \left( \prod_{i = 1}^{j-1} \fct{lm}_{\preceq}([\alpha_i , \alpha_j ]) \cdot \prod_{i = j+1}^n \fct{lm}_{\preceq}([ \alpha_i , \alpha_j ]) \right) \cdot \prod_{\substack{i \in I \\ k \in [n-3]} } A_i B_{k-1}^{-1}D_{k-1} \cdot \prod_{i \in I}A_i \nonumber \\
& = \prod_{\substack{ j \in J}} \left( \prod_{i = 1}^{j-1} A_iA_{j-1}^{-1}C_{j-1} \cdot \prod_{i = j+1}^n A_jA_{i-1}^{-1}C_{i-1} \right) \cdot \left( \prod_{ k \in [n-3]} B_{k-1}^{-1}D_{k-1}\right)^r \cdot \prod_{i \in I}A_i^{n-2}. \label{lmI}
}
Here, the first factor is
\eqns{
& \prod_{\substack{ j \in J}} \left( \prod_{i = 1}^{j-1} A_iA_{j-1}^{-1}C_{j-1} \cdot \prod_{i = j+1}^n A_jA_{i-1}^{-1}C_{i-1} \right) \\
& = \prod_{\substack{ j \in J}} (A_{j-1}^{-1}C_{j-1})^{j-1} A_j^{n-j} \left( \prod_{i = 1}^{j-1} A_i \cdot \prod_{i = j+1}^n A_{i-1}^{-1}C_{i-1} \right).}
The degree of $A_l$ in this factor is given by
\eqns{
&\deg_{A_l} \prod_{ j \in J} (A_{j-1}^{-1}C_{j-1})^{j-1} A_j^{n-j} \left( \prod_{i = 1}^{j-1} A_i \cdot \prod_{i = j+1}^n A_{i-1}^{-1}C_{i-1} \right) \\
& = \begin{cases} n-l & (l \in J) \\ 0 & (l \in I) \end{cases} - \begin{cases} l & (l+1 \in J) \\ 0 & (l+1 \in I) \end{cases} + \deg_{A_l}\prod_{ j \in J} \prod_{i = 1}^{j-1} A_i \cdot \prod_{i = j}^{n-1} A_i^{-1} \\
& = \begin{cases} n-l & (l \in J) \\ 0 & (l \in I) \end{cases} - \begin{cases} l & (l+1 \in J) \\ 0 & (l+1 \in I) \end{cases} + \begin{cases}\card{ ( J \cap [n] \setminus [l])} -  \card{ ( J \cap [l] ) } & (l \leq n-1) \\ 0 & (l = n)\end{cases}
}
Since we put $B_0 = A_n$, the degree of $A_l$ on the remaining factor is 
\eqns{
\deg_{A_l} \left( \prod_{ k \in [n-3]} B_{k-1}^{-1}D_{k-1}\right)^r \cdot \prod_{i \in I}A_i^{n-2} 
&= \begin{cases} -r & (l = n) \\ 0 & (l \leq n-1)\end{cases} + \begin{cases} n-2 & (l \in I) \\ 0 & (l \in J). \end{cases}
}
By summing up these degrees, we obtain the assertion.
\end{proof}
To apply Corollary \ref{mulindep} for discriminant-resultants, we use the expression of Theorem \ref{main1}. We firstly show that only the term of $I = [r]$ in the RHS of (\ref{eqmain1}) gives the maximal leading monomial in the terms of $DR_{n,r}$.
\begin{lemm}\label{degdiscres} For any integer $r=0$ or $2 \leq r \leq n$, we have
% \eqns{
% \deg_{A_l} \fct{lm}_{\preceq}(DR_{n,r}) & = \deg_{A_l} \fct{lm}_{\preceq} \left(  \prod_{\substack{ j \in [n] \setminus [r] \\ i \in [n] \setminus \{ j \}}} [\alpha_i , \alpha_j ] \cdot \prod_{\substack{i \in [r] \\ k \in [n-2]} } [\beta_k , \alpha_i] \right)\\
% & =  \begin{cases} 2n-r-2 & (1 \leq l \leq r-1) \\ 2n-2r-2 & (l = r) \\ 2n-4l+r & (r+1 \leq l \leq n-1) \\ 
% -r & (l = n,\ r \leq n-1).
% \end{cases}
% }
\eqn{
\fct{lm}_{\preceq}(DR_{n,r}) = \fct{lm}_{\preceq} \left(  \prod_{\substack{ j \in [n] \setminus [r] \\ i \in [n] \setminus \{ j \}}} [\alpha_i , \alpha_j ] \cdot \prod_{\substack{i \in [r] \\ k \in [n-2]} } [\beta_k , \alpha_i] \right).
}
\end{lemm}
\begin{proof} Since the monomial ordering $\preceq$ is the lexicographic order, it is enough to show the following claim for $1 \leq l \leq r-1$.
\begin{itemize}
    \item[($C_{r,l}$)] In the terms of the RHS of (\ref{eqmain1}) corresponding to the sets $I$ such that $I \supset [l-1]$, the degree of $A_l$ on their leading monomials is maximal if $l , l+1 \in I$.
\end{itemize}
By applying Lemma \ref{mondegAl} for the case $I \supset [l-1]$ and $r \leq l$, we obtain 
\eqns{
& \deg_{A_l} \fct{lm}_{\preceq} \left(  \prod_{\substack{ j \in J \\ i \in [n] \setminus \{ j \}}} [\alpha_i , \alpha_j ] \cdot \prod_{\substack{i \in I \\ k \in [n-2]} } [\beta_k , \alpha_i] \right)\\
& = 2n-r-2 - \begin{cases}0 & (l \in I) \\ l& (l \in J) \end{cases} - \begin{cases} 0 & (l+1 \in I) \\ l & (l+1 \in J). \end{cases}
}
This shows the claim $(C_{r,l})$.
\end{proof}
\begin{thm}(Theorem \ref{main3}) The discriminant-resultants
\[\{ DR_{n,r} \mid r=0,2,3,\ldots ,n \}\]
are algebraically independent. \label{algindep}
\end{thm}
\begin{proof} The assertion is obvious for $n = 2$, so we assume $n \geq 3$. By Corollary \ref{mulindep}, it is sufficient to show that the matrix \eqn{P := \left( \deg_{X} \fct{lm}_{\preceq}(DR_{n,r})\right)_{\substack{ r=0,2,\ldots ,n, \\ X: \text{variable of }L }}} is of full rank. By Lemma \ref{degdiscres} and (\ref{lmI}), 
\eqnsl{
& \fct{lm}_{\preceq} (DR_{n,r}) = \fct{lm}_{\preceq} \left(  \prod_{\substack{ j \in J \\ i \in [n] \setminus \{ j \}}} [\alpha_i , \alpha_j ] \cdot \prod_{\substack{i \in I \\ k \in [n-2]} } [\beta_k , \alpha_i] \right) \nonumber \\
& = \prod_{j = r+1}^n \left( \prod_{i = 1}^{j-1} A_iA_{j-1}^{-1}C_{j-1} \cdot \prod_{i = j+1}^n A_jA_{i-1}^{-1}C_{i-1} \right) \cdot \left( \prod_{ k = 1}^{n-3} B_{k-1}^{-1}D_{k-1}\right)^r \cdot \prod_{i = 1}^{r} A_i^{n-2}. \label{lmDR}
}
Using (\ref{lmDR}), we compute
\eqns{ c_{r,l} := \deg_{C_l} \fct{lm}_{\preceq} (DR_{n,r}) = \begin{cases} 0 & (l \leq r-1)\\
2l-r & (r \leq l \leq n-1) \\
r & (l = n)
\end{cases}
}
and
\eqns{ a'_{r,l} := & \deg_{A_l} \fct{lm}_{\preceq} (DR_{n,r}) - \deg_{C_l} \fct{lm}_{\preceq} (DR_{n,r}) \\
= &\deg_{A_l} \prod_{j = r+1}^n \left( \prod_{i = 1}^{j-1} A_i \cdot \prod_{i = j+1}^n A_j \right) \cdot \prod_{i = 1}^{r} A_i^{n-2}\\
= &
\begin{cases} 2n-r-2 & (l \leq r)\\
2n-2l & (r+1 \leq l \leq n).
\end{cases}
}
that the vectors
\[ v_r := ( a'_{r,2}- a'_{r,3} , c_{r,2}, \ldots , c_{r,n-1}) \]
for $r=0,2,3,\ldots , n$ are linear independent, because the matrix given by \eqn{(v_0 \  v_2 \ v_3 \  \cdots \  v_n)} is lower triagonal.
Thus, the matrix $P$ is of full rank. This is what we wanted to show.
% A = [[2 0][2 0]]

\end{proof}

\end{document}